\documentclass[a4paper, 11pt]{article}	 

\usepackage[paper=a4paper,left=1.8cm,right=1.8cm]{geometry}
\usepackage{cite}
\usepackage{colortbl}
\usepackage{color}
\usepackage{footnote}
\usepackage{subfig}
\usepackage{rotating} 
\usepackage[onehalfspacing]{setspace}
\usepackage{amsmath}
\usepackage{amsthm}
\usepackage{multirow}
\usepackage{amssymb}
\usepackage{psfrag}
\usepackage{graphicx}
\usepackage{algorithm}
\usepackage{algpseudocode}
\usepackage[T1]{fontenc}
\usepackage{authblk}
\usepackage[]{units}

\title{Simulation studies on regional predictive control}
\author{Kai K\"onig and Martin M\"onnigmann\thanks{Corresponding author.}\\
Automatic Control and Systems Theory, Department of Mechanical Engineering,\\
	Ruhr-Universit\"at Bochum, 44801 Bochum, Germany.\\ E-mail: {\tt\small kai.koenig-h4d@rub.de} and {\tt\small martin.moennigmann@rub.de}}

\newcommand{\R}{{\mathbb R}}

\newtheorem{lemma}{Lemma}

\begin{document}
\maketitle
\thispagestyle{plain}
\pagestyle{plain}
\section{Introduction}
We apply the regional predictive control approaches proposed in \cite{Jost2015a, Koenig2017a, Koenig2017c, Koenig2017d, Koenig2020, Moennigmann2020} to several examples to evaluate their efficiency. We first introduce the problem, the idea of regional predictive control and the existing approaches \cite{Jost2015a, Koenig2017a, Koenig2017c, Koenig2017d, Koenig2020, Moennigmann2020} to make this report self-contained. We then present numerical results as well as hardware-in-the-loop (HIL) results.  

\section{Brief problem statement}\label{sec:mpc}

We consider a linear discrete-time system
\begin{align}\label{eq:Sys}
x(k+1)&=Ax(k)+Bu(k), \ x(0) \ \text{given},
\end{align}
with state variables $x(k) \in \R^n$, input variables $u(k) \in \R^m$ and system matrices $A \in \R^{n \times n}$ and $B \in \R^{n \times m}$. We assume state and input constraints
\begin{align*}
x(k) \in \mathcal{X} \subset \R^n,\quad \quad
u(k) \in \mathcal{U} \subset \R^m
\end{align*}
apply for all $k \geq 0$. It is our objective to asymptotically stabilize the origin while satisfying the constraints. 
For this purpose, the OCP
\begin{align}\label{eq:MPCProblem}
\begin{split}
\min \limits_{X,U}  \quad &\tilde{x}(N)'P\tilde{x}(N)+\sum \limits_{i=0}^{N-1} (\tilde{x}(i)'Q\tilde{x}(i)+\tilde{u}(i)'R\tilde{u}(i))\\
\text{s.t.} \quad &\tilde{x}(0)= x,\\
&\tilde{x}(i+1)=A\tilde{x}(i)+B\tilde{u}(i), \quad i=0, \ldots, N-1,\\
&\tilde{x}(i) \in \mathcal{X}, \quad i=0, \ldots, N-1,\\
&\tilde{u}(i) \in \mathcal{U}, \quad i=0, \ldots, N-1,\\
&\tilde{x}(N) \in \mathcal{T}
\end{split}
\end{align}
is solved perpetually with the prediction horizon $N$, weighting matrices $Q \in \R^{n \times n}$, $R \in \R^{m \times m}$ and $P\in \R^{n \times n}$ and the terminal set $\mathcal{T}$. The solution of the OCP provides the optimal predicted state sequence $X=(\tilde{x}(1)', \ldots, \tilde{x}(N)')'$ and input sequence $U=(\tilde{u}(0)', \ldots, \tilde{u}(N-1)')'$ for the current system state $x=x(k)$. A closed-loop system results from applying the first $m$ elements of the input sequence, i.e. $\tilde{u}(0)$, to the system \eqref{eq:Sys}. We assume the matrices $Q$, $R$ and $P$ are symmetric and positive definite. Moreover, we assume the pair $(A,B)$ is stabilizable and the pair $(Q^{\frac{1}{2}}, A)$ is detectable. Finally, we assume $\mathcal{X}$, $\mathcal{U}$ and $\mathcal{T} \subseteq \mathcal{X}$ are convex and compact polytopes that contain the origin as an interior point. 
The weighting matrix $P$ and the terminal set $\mathcal{T}$ are chosen such that system \eqref{eq:Sys} is stabilized asymptotically while satisfying the constraints \cite{Mayne2000, Gilbert1991}. The matrix $P$ is the solution of the discrete-time algebraic Riccati equation. The terminal set $\mathcal{T}$ is calculated according to \cite{Gilbert1991}. By inserting the system dynamics \eqref{eq:Sys} into the cost function of the OCP \eqref{eq:MPCProblem}, the states $\tilde{x}(1), \dots, \tilde{x}(N)$ can be eliminated from the OCP. The resulting problem is a quadratic program (QP) of the form  
\begin{align}\label{eq:ReformulatedMPCProblem}
\begin{split}
\min \limits_U \ &\frac{1}{2}U'HU+x'FU+\frac{1}{2}x'Yx  \quad \\
&\text{s.t.} \quad GU \leq w+Ex
\end{split}
\end{align}
with $Y \in \R^{n \times n}$, $F \in \R^{n \times mN}$, $H \in \R^{mN \times mN}$, $G \in \R^{q \times mN}$, $w \in \R^q$, $E \in \R^{q \times n}$ and the number of constraints $q$. Note that $H\succ 0$, if $Q\succ 0$, $R\succ 0$ and $P\succ 0$. 
Let $\mathcal{X}_f$ refer to the set of initial states for which problem \eqref{eq:ReformulatedMPCProblem} has a solution. Under the assumptions stated for the problem \eqref{eq:MPCProblem}, $H$ is positive definite and there exists a unique optimal input sequence $U^\star(x)$ for every $x \in \mathcal{X}_f$. The optimal solution $U^\star: \mathcal{X}_f\rightarrow\R^{
mN}$ is a continuous piecewise affine function on a partition of $\mathcal{X}_f$ into a finite number of polytopes $\mathcal{P}^\star_1$, $\mathcal{P}^\star_2, \ldots$. We denote a single affine piece of the piecewise affine function by
\begin{align}\label{eq:FullAffineLaw} 
x \mapsto \bar{K}_j^\star x+ \bar{b}_j^\star \quad \forall x \in \mathcal{P}^\star_j 
\end{align}
with $\bar{K}^{\star}_j \in \R^{mN \times n}$ and $\bar{b}^{\star}_j \in \R^{mN}$, where we often omit the index $j$ for simplicity. Note that this affine function \eqref{eq:FullAffineLaw} yields the entire sequence of optimal signals that result from the optimal control problem. We refer to the first input signal, i.e., the first $m$ elements of \eqref{eq:FullAffineLaw}, as the feedback law for brevity, since it yields the MPC feedback signal on its polytope of validity. We denote the feedback law by $x \mapsto  K^\star x + b^\star \ \forall x \in \mathcal{P}^\star$ with $K^\star=\bar{K}^{\star {\{1, \ldots, m\}}}$ and $b^\star=\bar{b}^{\star {\{1, \ldots, m\}}}$, where a matrix and vector with a set index refer to the obvious submatrix and subvector. Proofs of the statements summarized in this section can be found in~\cite{Bemporad2002}.

\section{Regional predictive control}\label{sec:method}

\subsection{Basic approach \cite{Jost2015a}}

Regional predictive control makes use of the piecewise affine structure of the solution to the problem \eqref{eq:MPCProblem} (without calculating the solution explicitly for all feasible states). It exploits the fact that from the optimal solution of the OCP \textit{at a point} $x\in\mathcal{X}_f$ an optimal feedback law $K^\star x + b^\star$ and its polytope $\mathcal{P}^\star$ are known. This is stated more precisely in Lemma \ref{lem:Jost2015}, which bases on the results in \cite{Bemporad2002}. Lemma~\ref{lem:Jost2015} does not depend on the optimal solution $U^\star(x)$ for a point $x\in\mathcal{X}_f$ but only on the sets of active and inactive constraints 
\begin{align} \label{eq:Sets}
\begin{split}
\mathcal{A}(x)&=\{i \in \mathcal{Q}~|~G^iU^\star(x)=w^i+E^ix \},\\
\mathcal{I}(x)&=\{i \in \mathcal{Q}~|~G^iU^\star(x)<w^i+E^ix \}
\end{split}
\end{align}
with $\mathcal{Q}= \{1, \dots, q\}$ and $  \mathcal{I}(x)= \mathcal{Q}\backslash\mathcal{A}(x)$. Obviously, the sets $\mathcal{A}(x)$ and $\mathcal{I}(x)$, or $\mathcal{A}$ and $\mathcal{I}$ for short, can be determined by inserting the point $x\in\mathcal{X}_f$ and the solution $U^\star(x)$ in the constraints $G U\le w+ E x$ unless they are already available after solving \eqref{eq:MPCProblem} or \eqref{eq:ReformulatedMPCProblem}. We also need the weakly active set $\mathcal{W}(x)=\{i \in \mathcal{A}(x)~|~\lambda^{\star i}(x)=0 \}$ with the Lagrange multipliers $\lambda^\star \colon \mathcal{X}_f \rightarrow \R^q$ below.

\begin{lemma}\label{lem:Jost2015}
\cite{Jost2015a} Let $x \in \mathcal{X}_f$ be arbitrary and $\mathcal{A}(x)=\mathcal{A}$ the corresponding active set. Assume the matrix $G^{\mathcal{A}}$ has full row rank. Let 
\begin{align}\label{eq:LemJost}
  \begin{split}
    \bar{K}^\star&=H^{-1}(G^{\mathcal{A}})'(G^{\mathcal{A}}H^{-1}(G^{\mathcal{A}})')^{-1}S^{\mathcal{A}}-H^{-1}F',\\
    \bar{b}^\star&=H^{-1}(G^{\mathcal{A}})'(G^{\mathcal{A}}H^{-1}(G^{\mathcal{A}})')^{-1}w^{\mathcal{A}},\\
     T^\star&=\begin{pmatrix}G^{\mathcal{I}}H^{-1}(G^{\mathcal{A}})'(G^{\mathcal{A}}H^{-1}(G^{\mathcal{A}})')^{-1}S^{\mathcal{A}}-S^{\mathcal{I}} \\ (G^{\mathcal{A}}H^{-1}(G^{\mathcal{A}})')^{-1}S^{\mathcal{A}} \end{pmatrix},\\
      d^\star&=-\begin{pmatrix}G^{\mathcal{I}}H^{-1}(G^{\mathcal{A}})'(G^{\mathcal{A}}H^{-1}(G^{\mathcal{A}})')^{-1}w^{\mathcal{A}}-w^{\mathcal{I}} \\ (G^{\mathcal{A}}H^{-1}(G^{\mathcal{A}})')^{-1}w^{\mathcal{A}} \end{pmatrix},
  \end{split}
\end{align}
where  $S=E+GH^{-1}F'$, $S \in \R^{q \times n}$. Then the affine law $\bar{K}^\star x+\bar{b}^\star$ yields the optimal input sequence on the entire polytope 
$ \mathcal{P}^\star=\{ x \in \R^n ~|~T^\star x \leq d^\star \}$, i.e. $U^\star (x)=\bar{K}^\star x + \bar{b}^\star$ for all $x \in \mathcal{P}^\star$.
\end{lemma}

Lemma \ref{lem:Jost2015} suggest solving the OCP~\eqref{eq:MPCProblem} or QP~\eqref{eq:ReformulatedMPCProblem} only if the current feedback law loses its validity in the case of leaving the current polytope. For this, problem \eqref{eq:ReformulatedMPCProblem} is solved for the current state at first. With the current set of active constraints, the current optimal feedback law can be determined according to Lemma \ref{lem:Jost2015} and be reused until its polytope has been left. When leaving the polytope, the QP must be solved again and the feedback law and the polytope must be updated. If the rank condition $G^{{\mathcal{A}}}$ for a state $x$  is not met, the QP~\eqref{eq:ReformulatedMPCProblem} must be solved in the next time step.  

\subsection{New approaches}\label{sec:newApproaches}
It is the aim of the regional predictive control approach to reduce the number of optimization problems to be solved by reusing feedback laws whenever possible. We developed several approaches that increase the reusability further compared to the basic approach. We briefly introduce the new regional MPC approaches, which are compared to each other in the results sections below. We emphasize that all approaches are online approaches that \textit{do not} need the explicit solution to be known before. 
\paragraph*{{\bf Sets with common optimal feedback laws  \cite{Koenig2020EJC}}} Often a number of polytopes have the same optimal feedback law $K^\star x+b^\star$ in common. This is possible because a feedback law can be uniquely defined by a subset $\tilde{\mathcal{A}} \subseteq \mathcal{A}$ of the active set, which is the same on a union of polytopes. It is an obvious idea to reuse a feedback law not only on a single polytope $\mathcal{P}^\star$ as proposed in \cite{Jost2015a} 
but on this union of polytopes whenever possible. The maximum number of polytopes that are computed for a feedback law can be limited resulting in a heuristic approach (see \cite{Koenig2020EJC}, \S 4). We use the heuristic approach for the results presented below.    

\paragraph*{{\bf Active set updates \cite{Koenig2017a}}}
After leaving the current polytope, the new polytope and its feedback law can often be determined without solving a new QP. It is the idea of this approach to update the active set along a line connecting the current and the previous state. This is done by analyzing the crossed facets of neighboring polytopes along the line. By this, a number of neighboring polytopes and feedback laws along the closed-loop trajectory can be computed without solving a QP. A new QP has to be solved only if the linear independence constraint qualification (LICQ) is violated or more than one constraint is weakly active on a crossed facet.     

\paragraph*{{\bf Closed-loop optimal sequences of affine laws \cite{Monnigmann2019, Moennigmann2020}}}
With this approach all polytopes (and their feedback laws) that contain a state of the closed-loop trajectory can be computed from the solution of a QP at the current state, i.e., a single point $x \in \mathcal{X}_f$. If the terminal constraints are inactive at the current state $x$, then the solution of the QP for the state $x$ does not only provide a single feedback law but the entire sequence of optimal feedback laws and their polytopes of validity along the closed-loop trajectory. Consequently, a QP has to be solved only for states that result in active terminal constraints.

\paragraph*{{\bf Nonlinearly bounded regions of validity \cite{Koenig2017c}}}A feedback law $K^\star x+ b^\star$ can be used even if its polytope $\mathcal{P}^\star$ is left as long as it is feasible and stabilizing. With this insight, the validity of a feedback law can be extended from $\mathcal{P}^\star$ to a larger nonlinearly bounded region $\mathcal{E}$ that results by intersecting a polytopic feasibility region and a stability region defined by a simple quadric. The region $\mathcal{E}$ is computed instead of $\mathcal{P}^\star$ whenever possible in this approach. A QP is solved if the region of validity for a feedback law is left. 

\subsection{Realization in a networked setting}\label{sec:networkedSetting}
In \cite{BernerP2019d} the basic regional MPC approach \cite{Jost2015a} has been implemented in a networked MPC variant, where lean, low-power local hardware can be used to compute optimal closed-loop control signals. In this setting the QP \eqref{eq:ReformulatedMPCProblem} is solved on a computationally powerful central node on demand. The resulting active set according to \eqref{eq:Sets} is transmitted to a lean local node. On the local node, a feedback law and its polytopic region of validity are computed from the active set according to \eqref{eq:LemJost}. After that, the local node generates closed-loop input signals simply by evaluating the optimal affine feedback law. Whenever the region of validity for a feedback law is left, the central node is requested to solve a new QP. Note that transmitting active sets in a network is beneficial because they can be represented as tuples of $q$ bits $\alpha=(\alpha_q,\ldots,\alpha_1)$, where $\alpha_i=1$ if $i \in \mathcal{A}$ and $\alpha_i=0$ otherwise, leading to low bandwidth requirements.

We implement the new approaches from Section \ref{sec:newApproaches} in a networked MPC variant as in  \cite{BernerP2019d} in Section \ref{sec:results} below. In all approaches, QPs are solved on the central node on-demand and active sets are transmitted to the local node as in \cite{BernerP2019d}. Unlike the basic approach \cite{Jost2015a}, in the new approaches \cite{Koenig2020EJC, Koenig2017a, Moennigmann2020} the central node does not transmit a single set $\mathcal{A}_1$ but a number of sets $\mathcal{A}_1, \mathcal{A}_2, \ldots$, i.e., all sets that can be generated from the solution of a single QP (see \cite{Koenig2020EJC, Koenig2017a, Moennigmann2020} for details). The received active sets are used to compute feedback laws and polytopes as in the basic approach. Note that transmitting active sets is a good trade-off between the amount of transmitted data and the computational effort on the local node. For more details, we refer to \cite{BernerP2019d}.

\section{Numerical results} \label{sec:simstudy}
We apply the regional MPC approaches presented in the previous section to five numerical examples and compare the results. 

\subsection{Numerical examples}\label{subsec:examples}

\paragraph*{{\bf Artificial SISO (SISO20):}} 
	Consider the single-input-single-output system with the transfer function
	\begin{align}\label{example:SISO}
  	G(s)=\frac{2}{s^2+s+2} \nonumber
	\end{align}
	that is discretized with the 
	sampling time $T_s=\unit[0.1]{s}$. This results in a system of the form \eqref{eq:Sys} with 
	\begin{align}
		A=\begin{pmatrix}  
		0.8955 & -0.1897 \\
		0.0948 & 0.9903
		\end{pmatrix}, \quad
		B=\begin{pmatrix}
		0.0948 \\
		0.0048
		\end{pmatrix}. \nonumber
	\end{align}
	The example is similar to the one in~\cite{Seron2003}, but the system must here respect $-3 \leq x_i \leq 3$, $i=1,2$ 
	and $-2 \leq u_1 \leq 2$ 
	and weighting matrices read $Q=\text{diag}(0.01,4)$ and $R=0.01$. 
	We choose the horizon $N= 20$, which results in a QP with $q=128$ inequalities and 20 optimization variables.

\paragraph*{{\bf Ball and Plate (BP10):}} 
Consider a system of a ball rolling over a plate actuated by two independent motors. The state vector reads $x(k)=(y(k),\dot{y}(k),\alpha(k),\dot{\alpha}(k))'$ with ball position $y(k)$ and plate angle $\alpha(k)$. A discretization with the sampling time $T_s=\unit[0.03]{s}$ results in a system of the form \eqref{eq:Sys} with $n=4$, $m=1$ and matrices 
\begin{align}
A=\begin{pmatrix}
 1 & 0.0300 & 0.3150 & 0.0025\\
 0 & 1 & 21 & 0.2291\\
 0 & 0 & 1 & 0.0186\\
 0 & 0 & 0 & 0.3532
 \end{pmatrix}, \quad 
 B=\begin{pmatrix}
 0.00006\\
 0.0077 \\
 0.0010 \\
 0.0580 \\
 \end{pmatrix}. \nonumber
\end{align}
The state and input constraints read $-30 \leq y \leq 30$, $-15 \leq \dot{y} \leq 15$, $-0.26 \leq \alpha \leq 0.26$, $- 1 \leq \dot{\alpha} \leq 1$ and $-10 \leq u \leq 10$. The weighting matrices are set to $Q=\text{diag}(6,0.1,500,100)$ and $R=1$. We choose the horizon $N=10$, which results in a QP with $q=144$ inequalities and $10$ optimization variables. The system has been taken from \cite{Christophersen2007}.

\paragraph*{{\bf Inverted Pendulum (INPE50):}} 
Consider an inverted pendulum on a cart. The state vector reads $x(k)=(s(k),\varphi(k),\dot{s}(k),\dot{\varphi}(k))'$ with cart position $s(k)$ and pendulum angle $\varphi(k)$. A discretization with the sampling time $T_s=\unit[0.01]{s}$ results in a system of the form \eqref{eq:Sys} with $n=4$, $m=1$ and matrices 
\begin{align}
A=\begin{pmatrix}
 1 & -4.37 \cdot 10^{-5} & 0.0099 & 1.32 \cdot 10^{-7}\\
 0 & 1.0011 & 1.94 \cdot 10^{-4} & 0.0100\\
 0 & -0.0087 & 0.9812 &1.17 \cdot 10^{-5}\\
 0 & 0.2148 & 0.0386 & 0.9997
 \end{pmatrix}, \quad 
 B=\begin{pmatrix}
1.49 \cdot 10^{-5}\\
-3.08 \cdot 10^{-5}\\
 0.0030 \\
 -0.0061 \\
 \end{pmatrix}. \nonumber
\end{align}
The state and input constraints read $-1 \leq s \leq 1$, $- \frac{\pi}{3} \leq \varphi \leq \frac{\pi}{3}$, $-9 \leq \dot{s} \leq 9$, $- 2 \pi \leq \dot{\varphi} \leq 2 \pi$ and $-10 \leq u \leq 10$. The weighting matrices are set to $Q=I^{4 \times 4}$ and $R=0.01$. We choose the horizon $N=50$, which results in a QP with $q=804$ inequalities and $50$ optimization variables.

\paragraph*{{\bf Connected masses (COMA40):}} This system contains a chain of six masses connected to each other by springs, and to rigid walls on both ends of the chain. All masses and all spring constants are set to unity. Three forces $u_1, u_2, u_3$ between the first and second, third and fifth, and fourth and sixth mass, respectively, are used as inputs. 
The resulting system has $12$ states and $3$ inputs. Discretizing with a sampling time $T_s=\unit[0.5]{s}$ yields a system of the form \eqref{eq:Sys}. The system matrices are given in the appendix. The state and input constraints read $-4\leq x_i(t)\leq 4, \, i= 1, \dots, 12$ and $-0.5 \leq u_j(t)\leq 0.5, \, j= 1, \dots, 3$, respectively. We choose $Q = I^{12\times 12}$ and $R = I^{3\times 3}$.
The resulting QP has $mN= 120$ decision variables and $q=1282$ inequality constraints for a horizon of $N= 40$. The system has been taken from \cite{Boyd2010}.

\paragraph*{{\bf Artificial MIMO (MIMO75):}} Consider the multiple-input-multiple-output system with the transfer function
 \begin{align}
  G(s) = 
  \begin{pmatrix}
   \frac{-5s+1}{36s^2+6s+1} & \frac{0.5s}{8s+1}           & 0\\
   0                        & \frac{0.1(-10s+1)}{s(8s+1)} & \frac{-0.1}{(64s^2+6s+1)s}\\
   \frac{-2s+1}{12s^2+3s+1} &  0& \frac{2(-5s+1)}{16s^2+2s+1}
  \end{pmatrix} \nonumber
 \end{align}
 that is discretized with the sampling time $T_s=\unit[1]{s}$. After removing uncontrollable states a system of the form \eqref{eq:Sys} with $n= 10$ states and $m= 3$ inputs results. 
The state and input constraints read $-10\leq \, x_i(t) \leq 10, \, i= 1, \dots, 10$ and $-1\leq \,u_j(t) \leq 1, \, j= 1, \dots, 3$ for this example. Furthermore, $Q = I^{n\times n}$, $R = 0.25I^{m\times m}$ and $N = 75$. The resulting QP has $mN= 225$ decision variables and $q= 2092$ inequality constraints. The system is taken from \cite{jost2014simulation}. 

\subsection{Numerical results}
We compare the regional MPC approaches presented in the previous section in terms of the reusability of the feedback laws and the computational effort. For this, we generate 10000 random initial states $x \in \mathcal{X}_f$ for every system and compute trajectories for the MPC-controlled system until $||x(k)|| \leq 10^{-3}$. For a given system, we use the same initial conditions for all approaches. Table \ref{tab:resultsAP1u4} shows the results for the sample systems from Section \ref{subsec:examples}. Here, the reusability refers to the fraction of time steps up to the terminal set, in which a feedback law from the previous time step can be reused and thus the solution of a QP can be avoided, in percent. A value of \unit[100]{\%} means that, after solving a QP for the initial state, a feedback law from the previous time step can be reused in all time steps along the closed-loop trajectory. The computational effort was determined by measuring \textsc{matlab} execution times. It shows the reduction in relation to the basic approach in percent. The last column shows the average values across all sample systems. 

The results show almost no feedback law can be reused (\unit[0.3]{\%} on average) with the basic approach \cite{Jost2015a} (second row). In contrast, a significant increase of the reusability and thus a reduction in the computational effort may be achieved with the new approaches. 

Exploiting sets with common optimal feedback laws (third row) results in an average reusability of \unit[20.6]{\%}. The computational effort is reduced by \unit[8.4]{\%} on average. The reusability varies from \unit[0]{\%} for COMA40 and MIMO75 to \unit[45.6]{\%} for INPE50. The computational effort varies from \unit[-32.2]{\%} for INPE50 to even \unit[+3.9]{\%} for COMA40. It is interesting that the two largest systems have the worst results, but the best results are achieved by the third largest system. This shows that the effectiveness of this approach strongly depends on the specific system. Note that we limited the number of sets that are computed for a feedback law to keep the computational effort reasonable (see \cite{Koenig2020EJC} for details). 

A strong dependency on the specific system also exists for the approach using nonlinearly bounded regions of validity (last row). The reusability varies from \unit[5.8]{\%} for COMA40 to \unit[88.8]{\%} for INPE50. However, the computational effort can be reduced for all systems and varies from \unit[-4.4]{\%} for COMA40 to \unit[-72.5]{\%} for INPE50. The average reusability is increased by \unit[43.9]{\%} resulting in a reduction of the computational effort by \unit[33.1]{\%}. 

The best results are achieved with active set updates and closed-loop optimal sequences of affine laws.

If active set updates are used (fourth row) an average reusability of \unit[96.4]{\%} can be achieved. This reduces the computational effort by \unit[66.2]{\%}. Moreover, the approach shows good results for all systems. 
The same applies to the approach using closed-loop optimal sequences (fifth row). Here, average reusability of \unit[91.8]{\%} can be achieved. Although this is a little less than in the active set update approach, the reduction in the computational effort is larger with \unit[78.8]{\%} on average. For COMA40, the computational effort can even be reduced by \unit[91.7]{\%}. 

\begin{table}[b]
\scriptsize
\setlength{\tabcolsep}{1mm}
\centering
\begin{tabular}
{>{\centering }m{3.2cm} | >{\centering}m{0.5cm} | >{\centering}m{1.7cm} ||>{\centering}m{1.5cm} |>{\centering}m{1.5cm} | >{\centering}m{1.5cm} | >{\centering}m{1.5cm} |  >{\centering}m{1.5cm} || >{\centering}m{1.5cm}}

 approaches& 
ref.& 
measurand&
SISO20 &
BP10 &
INPE50 &
COMA40 &
MIMO75 &
\textbf{average}  \tabularnewline \hline \hline
basic approach & \cite{Jost2015a} & reusability &\unit[0.2]{\%}& \unit[0]{\%}&  \unit[1.2]{\%}& \unit[0]{\%} & \unit[0]{\%}& \unit[0.3]{\%}\tabularnewline \hline  
\multirow{2}{3cm}{sets with common optimal feedback laws} & \multirow{2}{*}{\cite{Koenig2020EJC}} &reusability& \unit[26.4]{\%}& \unit[30.9]{\%}&  \unit[45.6]{\%}&\unit[0]{\%}&\unit[0]{\%}&\unit[20.6]{\%}\tabularnewline 
 & & comp. effort & \unit[-5.7]{\%}&  \unit[-9.6]{\%}&\unit[-32.2]{\%}&\unit[3.9]{\%}&\unit[1.7]{\%}&\unit[-8.4]{\%}\tabularnewline \hline 
\multirow{2}{3cm}{active set updates} & \multirow{2}{*}{\cite{Koenig2017a}} & reusability & \unit[90.4]{\%} &\unit[99]{\%}& \unit[98.3]{\%}& \unit[97.4]{\%} &\unit[97.1]{\%}& \unit[96.4]{\%}  \tabularnewline 
 & & comp. effort & \unit[-50.5]{\%}&  \unit[-58.7]{\%}&\unit[-84.3]{\%}&\unit[-51.4]{\%}&\unit[-86]{\%}&\unit[-66.2]{\%}\tabularnewline \hline 
\multirow{2}{3.2cm}{closed-loop optimal sequences of affine laws} &\multirow{2}{*}{\cite{Moennigmann2020}}&reusability & \unit[94.9]{\%}&\unit[72.4]{\%}&\unit[97.9]{\%} &\unit[96.7]{\%}&\unit[96.9]{\%}&\unit[91.8]{\%}\tabularnewline 
 & & comp. effort & \unit[-74.5]{\%}&  \unit[-54.6]{\%}&\unit[-84.9]{\%}&\unit[-91.7]{\%}&\unit[-88.5]{\%}&\unit[-78.8]{\%}\tabularnewline \hline 
\multirow{2}{3cm}{nonlinearly bounded regions of validity} &\multirow{2}{*}{\cite{Koenig2017c}}&reusability & \unit[38.2]{\%}&\unit[44.9]{\%}&\unit[88.8]{\%} &\unit[5.8]{\%}&\unit[41.9]{\%}&\unit[43.9]{\%}\tabularnewline 
 & & comp. effort & \unit[-33.5]{\%}&  \unit[-27.3]{\%}&\unit[-72.5]{\%}&\unit[-4.4]{\%}&\unit[-27.9]{\%}&\unit[-33.1]{\%}\tabularnewline \hline 
 \end{tabular}
\caption{Numerical results. }\label{tab:resultsAP1u4}
\end{table}

\section{Hardware-in-the-loop results}\label{sec:results}
All results reported so far were obtained in simulations. We report results for an implementation on embedded hardware in this section. We use the networked setting described in Section \ref{sec:networkedSetting}. We use a standard desktop computer containing an Intel Core2 Duo CPU with two 3.0 GHz cores and 8 GB RAM as the central node. The central node is connected to an IEEE 802.11 b/g/n wireless LAN access point. We use Espressif ESP8266 SoC with an integrated IEEE 802.11 b/g/n WiFi controller as the local node. The SoC features an 80MHz Tensilica L106 32-bit RISC microcontroller and 96 KiB data RAM. It is connected to a dSpace hardware-in-the-loop simulator, which emulates the sample systems. We measure the states after analog to digital conversion and generate control inputs by digital to analog converters. In both cases, the resolution is 12 bit. We used the same experimental setup as explained in \cite[\S 4.3]{BernerDiss2019}. For more technical details we refer to the explanations therein. 

\subsection{HIL examples}\label{subsec:examplesHIL}
We compare the regional MPC approaches in a networked setting. For the comparison, we choose three examples with a lower problem size compared to the examples from Section \ref{sec:simstudy}. The reason for this is the limited memory of the ESP8266.   

\paragraph*{{\bf Double Integrator (DI6):}} 
	Consider the double integrator system that results from discretizing the continuous-time system
	\begin{align}
		\dot{x}_c(t)=\begin{pmatrix}  
		-1 & -2 \\
		1 & 0
		\end{pmatrix} x_c(t)+
		\begin{pmatrix}
		1 \\
		0
		\end{pmatrix} u_c(t) \nonumber
	\end{align}
	with the sampling time $T_s=\unit[1]{s}$. The system must respect $-3 \leq x_i \leq 3$, $i=1,2$ 
	and $-2 \leq u_1 \leq 2$ 
	and weighting matrices read $Q=\text{diag}(0.01,4)$ and $R=0.01$. 
	We choose the horizon $N= 6$, which results in a QP with $q=42$ inequalities and 6 optimization variables.
	
\paragraph*{{\bf Unstable System (US12):}} 
	Consider the unstable system that results from discretizing the continuous-time system
	\begin{align}
		\dot{x}_c(t)=\begin{pmatrix}  
		-1 & 0.3 \\
		0.1 & 1
		\end{pmatrix} x_c(t)+
		\begin{pmatrix}
		0.5 \\
		-2
		\end{pmatrix} u_c(t) \nonumber
	\end{align}
	with the sampling time $T_s=\unit[0.5]{s}$. The system must respect $-3 \leq x_i \leq 3$, $i=1,2$ 
	and $-1 \leq u_1 \leq 1$ 
	and weighting matrices read $Q=I^{2 \times 2}$ and $R=0.1$. 
	We choose the horizon $N= 12$, which results in a QP with $q=76$ inequalities and 12 optimization variables.
	
\paragraph*{{\bf Artificial MIMO (AM4):}} 
	Consider an artificial system with multiple inputs and outputs. The system results from discretizing the transfer function
	 \begin{align}
  G(s) = 
  \begin{pmatrix}
   \frac{0.5}{36s^2+6s+1} & \frac{0.02(2s+1)}{8s+1}     \\
   \frac{0.02(2s+1)}{8s+1} & \frac{0.05}{12s^2+3s+1}
  \end{pmatrix}, \nonumber
 \end{align}
	with the sampling time $T_s=\unit[1]{s}$. The system must respect $-10 \leq x_i \leq 10$, $i=1,2, \ldots, 6$ 
	and $-1 \leq u_j \leq 1$, $i=1,2$ 
	and weighting matrices read $Q=I^{6 \times 6}$ and $R=0.25 I^{2 \times 2}$. 
	We choose the horizon $N= 4$, which results in a QP with $q=104$ inequalities and 8 optimization variables.		
	
\subsection{Hardware-in-the-loop results}
We compare the new approaches to the basic approach \cite{Jost2015a} in terms of requests to the central node, amount of transmitted data, overall computational effort and the computational effort on the local node. For this, we generate closed-loop trajectories for 1000 random initial states for each sample system from Section \ref{subsec:examplesHIL}. We transmit active sets as bit tuples across the network as described in Section \ref{sec:networkedSetting}. The amount of transmitted data (rows labeled {\it{data}} in Table \ref{tab:resultsAP2u3u7}) is measured by recording the bits from the transmitted active sets. The computational effort is determined by measuring computation times on the microcontroller. {\it{Local effort}} refers to the time for computing and evaluating feedback laws and their regions of validity from the received active sets. {\it{Overall effort}} includes the time in which the local node is waiting for active sets from the central node. Results are recorded until the system enters the terminal region. The recorded values are the average values related to a single closed-loop trajectory. The last column of Table \ref{tab:resultsAP2u3u7} shows the average values across all sample systems.

In the basic approach \cite{Jost2015a} the average number of requests to the central node is 3.537 per trajectory. Along the trajectory, only 36.4593 bytes have to be transmitted to the local node on average. 

Exploiting sets with common optimal feedback laws reduces the number of requests by \unit[13.5]{\%}. Since not a single but a number of polytopes are transmitted to the local node per request, the amount of transmitted data is increased by \unit[65.2]{\%} and the computational effort on the local node is increased by \unit[28.0]{\%}. The overall effort can be reduced only for US12 with \unit[-24.9]{\%}. Again we limited the number of sets that are computed for a feedback law to keep the computational effort reasonable (see \cite{Koenig2020EJC} for details).   

The active set update approach may also transmit more than one polytope per request. The amount of transmitted data increase by \unit[34.6]{\%} and \unit[25.5]{\%}, respectively. However, the number of requests reduce by \unit[68.8]{\%} and \unit[52.6]{\%}, respectively. 

With closed-loop optimal sequences of affine laws, the number of requests reduce by \unit[65.3]{\%} and \unit[52.9]{\%}, respectively. The amount of transmitted data only slightly increase by \unit[3.8]{\%} and the local effort increase by \unit[11.2]{\%}.

Using nonlinearly bounded regions reduces the number of requests and the amount of transmitted data by \unit[16.9]{\%}. A reduction in the number of requests results in a reduction of the transmitted data because only a single active set is transmitted per request. However, the reduction in the number of requests is lower than in the approaches \cite{Koenig2017a, Moennigmann2020}. The local effort is increased by \unit[67.3]{\%} and the overall effort is reduced only by \unit[0.8]{\%}. The increased computational effort on the local node is related to the computation of the nonlinearly bounded regions (see \cite{Koenig2017c} for more details). 

The results show that the regional MPC approaches are appropriate for use in a networked setting with embedded hardware. Feedback laws and their regions of validity can be transmitted efficiently across the network by using an active set representation. The most suitable approach is based on the respective goal. Note that computations can systematically be shifted from the local to the central node and vice versa to affect the amount of transmitted data and the computation time on the local node (see \cite{BernerP2019d} for more details).        

\begin{table}[t]
\scriptsize
\setlength{\tabcolsep}{1mm}
\centering
\begin{tabular}
{>{\centering }m{2.3cm} | >{\centering}m{1.2cm} | >{\centering}m{1.8cm} ||>{\centering}m{1.7cm} | >{\centering}m{1.7cm} |  >{\centering}m{1.7cm} || >{\centering}m{1.7cm}}

 approach& 
 references&
measurand & 
DI6 &
US12 &
AM4 &
\textbf{average}  \tabularnewline \hline \hline
\multirow{2}{*}{basic approach} &\multirow{2}{*}{\cite{Jost2015a} }&requests& 2.663  & 3.306 & 4.642& 3.537 \tabularnewline \cline{3-7}
& & data [bytes] & 15.972  & 33.06 &60.346&36.4593\tabularnewline \hline  \hline 
\multirow{4}{2.3cm}{\centering sets with common optimal feedback laws} &\multirow{4}{*}{\cite{Koenig2020EJC}} &requests & \unit[-0]{\%}&{\unit[-33.6]{\%}}&  {\unit[-6.9]{\%}}&{\unit[-13.5]{\%}} \tabularnewline \cline{3-7}
 && data &  \unit[+9.3]{\%}&\unit[+73.8]{\%} & \unit[+112.4]{\%}&\unit[+65.2]{\%}\tabularnewline\cline{3-7}
   && overall effort  & \unit[+3.3]{\%} & {\unit[-24.9]{\%}}&\unit[+3.3]{\%} &   {\unit[-6.1]{\%}} \tabularnewline \cline{3-7}
 && local effort & \unit[+20.9]{\%} & \unit[+28.5]{\%}& \unit[+34.5]{\%}&   \unit[+28.0]{\%} \tabularnewline \hline 
 \multirow{4}{2.3cm}{\centering active set updates} &\multirow{4}{*}{\cite{Koenig2017a}}&requests&  {\unit[-62.0]{\%}} & {\unit[-69.8]{\%}}  &{\unit[-74.7]{\%}} &  {\unit[-68.8]{\%}} \tabularnewline \cline{3-7}
&& data & \unit[+40.8]{\%}  &\unit[+2.8]{\%}   & \unit[+60.1]{\%}&\unit[+34.6]{\%}  \tabularnewline\cline{3-7}
  & & overall effort &{\unit[-55.4]{\%}} & {\unit[-58.7]{\%}}&{\unit[-43.7]{\%}}& {\unit[-52.6]{\%}}     \tabularnewline \cline{3-7}
 &&local effort&\unit[+7.3]{\%}& \unit[+5]{\%}  &\unit[+64.3]{\%}&\unit[+25.5]{\%}     \tabularnewline \hline 
 \multirow{4}{2.3cm}{\centering closed-loop optimal sequences of affine laws} &\multirow{4}{*}{\cite{Moennigmann2020}}&requests&  {\unit[-61.7]{\%}} & {\unit[-69.8]{\%}}  &{\unit[-64.3]{\%}} &  {\unit[-65.3]{\%}} \tabularnewline \cline{3-7}
&& data & \unit[+1.5]{\%}  &\unit[+2.9]{\%}   & \unit[+6.9]{\%}&\unit[+3.8]{\%}  \tabularnewline\cline{3-7}
  & & overall effort & {\unit[-56.4]{\%}} & {\unit[-58.9]{\%}}&{\unit[-43.4]{\%}}& {\unit[-52.9]{\%}}     \tabularnewline \cline{3-7}
 && local effort&\unit[-1.82]{\%}& \unit[+5.4]{\%}  &\unit[+30.0]{\%}&\unit[+11.2]{\%}     \tabularnewline \hline 
  \multirow{4}{2.3cm}{\centering nonlinearly bounded regions} &\multirow{4}{*}{\cite{Koenig2017c}}&requests& \unit[+0.3]{\%} & {\unit[-15.1]{\%}}& {\unit[-35.9]{\%}}& {\unit[-16.9]{\%}} \tabularnewline \cline{3-7}
&& data& \unit[+0.3]{\%}  & {\unit[-15.1]{\%}}&{\unit[-35.9]{\%}}&{\unit[-16.9]{\%}}\tabularnewline\cline{3-7}
  && overall effort & \unit[+4.5]{\%} & {\unit[-10.7]{\%}}&\unit[+3.9]{\%}&    {\unit[-0.8]{\%}} \tabularnewline \cline{3-7}
 && local effort &  \unit[+38.2]{\%}  & \unit[+14.5]{\%}& \unit[+149.3]{\%} &    \unit[+67.3]{\%}  \tabularnewline \hline  \hline 
 \end{tabular}
\caption{Hardware-in-the-loop results.}  \label{tab:resultsAP2u3u7}
\end{table}  

\section{Conclusion}\label{sec:conclusion}
Regional predictive control as in \cite{Jost2015a} aims at reducing the number of optimization problems to be solved by reusing optimal affine feedback laws whenever possible. We applied the regional predictive control approaches proposed in \cite{Jost2015a, Koenig2017a, Koenig2017c, Koenig2017d, Koenig2020, Moennigmann2020} to several examples to evaluate their efficiency. We showed that the new approaches perform better compared to the basic regional MPC approach \cite{Jost2015a}. More precisely, the reusability of a feedback law can be increased and, thus, the number of QPs that must be solved can be reduced. The approaches based on active set updates \cite{Koenig2017a} and closed-loop optimal sequences of affine laws \cite{Moennigmann2020} perform best. The reusability and, thus, the reduction in the number of QPs is very high for all considered systems ($ > \unit[90]{\%}$ on average) even for large systems with more than 2000 constraints and more than 200 optimization variables. A significant reduction of the computational effort is achieved as a result.

Moreover, we showed that all approaches are appropriate for the implementation on lean embedded hardware in a networked MPC setting. Despite memory limitations all MPC approaches could be implemented on a lean ESP8266 microcontroller. The approach based on closed-loop optimal sequences performs best. It reduces the number of requests to the central node and the overall effort significantly while the computational effort on the local node and the amount of transmitted data is only slightly increased.    

\newpage
\section*{Appendix} 
The system matrices for the connected masses (COMA40) example are
 \begin{align*}
 A=[A_1, A_2]
 \end{align*}
 with
{\tiny
 \begin{align*}
A_1=[0.762721047593857,0.114882546593898,0.002476544740668,0.000020938074941,0.000000094222942,0.000000000263060, \\  0.114882546593898,0.765197592334525,0.114903484668839,0.002476638963610,0.000020938338001,0.000000094222942,\\  0.002476544740668,0.114903484668839,0.765197686557467,0.114903484931899,0.002476638963610,0.000020938074941,\\  0.000020938074941,0.002476638963610,0.114903484931899,0.765197686557467,0.114903484668839,0.002476544740668,\\  0.000000094222942,0.000020938338001,0.002476638963610,0.114903484668839,0.765197592334525,0.114882546593898,\\  0.000000000263060,0.000000094222942,0.000020938074941,0.002476544740668,0.114882546593898,0.762721047593857,\\                                                                                                             \text{-}0.899414767742325,0.420238976357871,0.019312099175890,0.000248250167124,0.000001497064606,0.000000005237232,\\                                                                                                   0.420238976357871,\text{-}0.880102668566434,0.420487226524995,0.019313596240496,0.000248255404355,0.000001497064606,\\                                                                                                      0.019312099175890,0.420487226524995,\text{-}0.880101171501828,0.420487231762227,0.019313596240496,0.000248250167124,\\  0.000248250167124,0.019313596240496,0.420487231762227,\text{-}0.880101171501828,0.420487226524995,0.019312099175890,\\  0.000001497064606,0.000248255404355,0.019313596240496,0.420487226524995,\text{-}0.880102668566434,0.420238976357871,\\  0.000000005237232,0.000001497064606,0.000248250167124,0.019312099175890,0.420238976357871,\text{-}0.899414767742325
]
\end{align*}
}
{\tiny
 \begin{align*}
A_2=[0.459613939727602,0.019813111712880,0.000251260056029,0.000001507575068,0.000000005261230,0.000000000011999,\\  0.019813111712880,0.459865199783630,0.019814619287947,0.000251265317259,0.000001507587067,0.000000005261230,\\  0.000251260056029,0.019814619287947,0.459865205044861,0.019814619299947,0.000251265317259,0.000001507575068,\\  0.000001507575068,0.000251265317259,0.019814619299947,0.459865205044861,0.019814619287947,0.000251260056029,\\  0.000000005261230,0.000001507587067,0.000251265317259,0.019814619287947,0.459865199783631,0.019813111712880,\\  0.000000000011999,0.000000005261230,0.000001507575068,0.000251260056029,0.019813111712880,0.459613939727602,\\                                                                                                0.762721047593857,0.114882546593898,0.002476544740668,0.000020938074941,0.000000094222942,0.000000000263060,\\                                                                        0.114882546593898,0.765197592334525,0.114903484668839,0.002476638963610,0.000020938338001,0.000000094222942,\\                                                        0.002476544740668,0.114903484668839,0.765197686557467,0.114903484931899,0.002476638963610,0.000020938074941,\\  0.000020938074941,0.002476638963610,0.114903484931899,0.765197686557467,0.114903484668839,0.002476544740668,\\  0.000000094222942,0.000020938338001,0.002476638963610,0.114903484668839,0.765197592334525,0.114882546593898,\\  0.000000000263060,0.000000094222942,0.000020938074941,0.002476544740668,0.114882546593898,0.762721047593857
]
\end{align*}
}
and
{\tiny
\begin{align*}
B = [ &0.117380123896010,   0.000021127047947,   0.000000094750065, \\
  &\text{-}0.117401251208021,   0.002518704613621,   0.000021127312010, \\
  &\text{-}0.002497672052678,   0.119898828510133,   0.002518704614624, \\
  &\text{-}0.000021032825507,   0.000000000000501,   0.119898828774196, \\
  &\text{-}0.000000094487005,  \text{-}0.119898828246070,   0.000000094750566, \\
  &\text{-}0.000000000263561,  \text{-}0.002518609863556,  \text{-}0.119877701198123, \\
   &0.439800828014722,   0.000251254794798,   0.000001507563068, \\
  &\text{-}0.440052088070751,   0.019813111700880,   0.000251260056029, \\
  &\text{-}0.019563359231919,   0.459613939727602,   0.019813111724879, \\
  &\text{-}0.000249757742191,   0.000000000011999,   0.459613944988832, \\
  &\text{-}0.000001502325837,  \text{-}0.459613934466372,   0.000001507575068, \\
  &\text{-}0.000000005249231,  \text{-}0.019811604137812,  \text{-}0.459362679671573] 
\end{align*}
}  

\section*{Acknowledgement}
Support by the Deutsche Forschungsgemeinschaft (DFG) under grant MO 1086/15-1 is gratefully acknowledged.

\bibliographystyle{plain}
\bibliography{literatur}

\begin{thebibliography}{10}

\bibitem{Bemporad2002}
A.~Bemporad, M.~Morari, V.~Dua, and E.~N. Pistikopoulos.
\newblock The explicit linear quadratic regulator for constrained systems.
\newblock {\em Automatica}, 38:3--20, 2002.

\bibitem{BernerDiss2019}
P.~Berner.
\newblock {\em An Event-Triggered Networked Model Predictive Control Approach
  for Lean Embedded Hardware}.
\newblock PhD thesis, 2019.

\bibitem{BernerP2019d}
P.~Berner and M.~Mönnigmann.
\newblock A complexity analysis of event-triggered model predictive control on
  industrial hardware.
\newblock {\em IEEE Transactions on Control Systems Technology}, pages 1--8,
  2019.

\bibitem{Christophersen2007}
F.~J. Christophersen, M.~Kvasnica, C.N. Jones, and M.~Morari.
\newblock {\em Optimal control of constrained piecewise affine systems}, volume
  359, chapter Efficient Evaluation of Piecewise Control Laws Defined Over a
  Large Number of Polyhedra, pages 150--165.
\newblock Springer, Berlin, 2007.

\bibitem{Gilbert1991}
E.~G. Gilbert and K.~T. Tan.
\newblock Linear systems with state and control constraints: the theory and
  application of maximal output admissible sets.
\newblock {\em IEEE Transactions on Automatic Control}, 36(9):1008--1020, 1991.

\bibitem{jost2014simulation}
M.~Jost, G.~Pannocchia, and M.~Mönnigmann.
\newblock Simulation studies on online constraint removal with a {L}yapunov
  function, 2014.

\bibitem{Jost2015a}
M.~Jost, M.~{Schulze Darup}, and M.~Mönnigmann.
\newblock Optimal and suboptimal event-triggering in linear model predictive
  control.
\newblock In {\em Proceedings of the 2015 European Control Conference}, pages
  1147--1152, 2015a.

\bibitem{Koenig2020}
K.~König and M.~Mönnigmann.
\newblock Accelerating {MPC} by online detection of state space sets with
  common optimal feedback laws.
\newblock {\em Optimal Control Applications and Methods (submitted)}.

\bibitem{Koenig2017a}
K.~König and M.~Mönnigmann.
\newblock Regional {MPC} with active set updates.
\newblock {\em IFAC-PapersOnLine}, 50(1):11859--11864, 2017.

\bibitem{Koenig2017d}
K.~König and M.~Mönnigmann.
\newblock {R}egionale prädiktive {R}egelung - {M}odellprädiktive {R}egelung
  mittels stückweise definiertem {R}iccati-{R}egler.
\newblock {\em at-Automatisierungstechnik}, 65(12):808--821, 2017.

\bibitem{Koenig2017c}
K.~König and M.~Mönnigmann.
\newblock Regional {MPC} with nonlinearly bounded regions of validity.
\newblock In {\em Proceedings of the 2018 European Control Conference (ECC)},
  pages 294--299, 2018.

\bibitem{Koenig2020EJC}
K.~König and M.~Mönnigmann.
\newblock Accelerating {MPC} by online detection of state space sets with
  common optimal feedback laws.
\newblock {\em Optimal Control Applications and Methods (submitted)}, 2020.

\bibitem{Mayne2000}
D.~Q. Mayne, J.~B. Rawlings, C.~V. Rao, and P.~O.~M. Scokaert.
\newblock Constrained model predictive control: Stability and optimality.
\newblock {\em Automatica}, 36:789--814, 2000.

\bibitem{Monnigmann2019}
M.~Mönnigmann.
\newblock On the structure of the set of active sets in constrained linear
  quadratic regulation.
\newblock {\em Automatica}, {106}:{61--69}, {2019}.

\bibitem{Moennigmann2020}
M.~Mönnigmann and G.~Pannocchia.
\newblock Reducing the computational effort of {MPC} with closed-loop optimal
  sequences of affine laws.
\newblock In {\em Proceedings of the 21st IFAC World Congress}, pages
  11508--11513, 2020.

\bibitem{Seron2003}
M.~M. Seron, G.~C. Goodwin, and J.~A. {De Doná}.
\newblock Characterisation of receding horizon control for constrained linear
  systems.
\newblock {\em Asian Journal of Control}, 5(2):271--286, 2003.

\bibitem{Boyd2010}
Y.~Wang and S.~Boyd.
\newblock Fast model predictive control using online optimization.
\newblock {\em IEEE Transactions on Control Systems Technology}, 18:267--278,
  2010.

\end{thebibliography}

\end{document}